\documentclass[pdftex]{imsart}

\RequirePackage[OT1]{fontenc}
\usepackage{a4wide}

\usepackage{hyperref}

\usepackage{amsmath,amsthm}
\usepackage
{inputenc}
\usepackage{amsfonts}
\usepackage{graphicx}
\usepackage{natbib}
\usepackage{amsthm}
\newtheorem{theorem}{Theorem}[section]

\newtheorem{corollary}[theorem]{Corollary}
\newtheorem{lemma}[theorem]{Lemma}

\theoremstyle{remark}
\newtheorem{remark}[theorem]{Remark}
\newtheorem{example}[theorem]{Example}
\theoremstyle{definition}

\startlocaldefs
\usepackage{url}		
\usepackage{xspace}		
\usepackage{booktabs}

\usepackage{rotating}	
\usepackage{pdflscape}
\usepackage{fancyvrb}

\usepackage{bbm}			

\newcommand{\SU}{\mathbf{SU}}
\newcommand{\SD}{\mathbf{SD}}

\newcommand{\DBR}{\textnormal{DBR}}

\newcommand{\RBH}{\textnormal{RBH}}

\newcommand{\FDR}{\textnormal{FDR}}

\newcommand{\Fbar}{\overline{F}}
\newcommand{\FbarSD}{\overline{F}_{\textnormal{SD}}}
\newcommand{\FbarSU}{\overline{F}_{\textnormal{SU}}}

\newcommand{\wt}{\widetilde}
\newcommand{\wh}{\widehat}
\newcommand{\mbf}{\mathbf}
\newcommand{\E}{\mathbf{E}}
\renewcommand{\P}{\mathbf{P}}
\newcommand{\cH}{\mathcal{H}}
\newcommand{\ind}[1]{{\mbf{1}\{#1\}}}

\endlocaldefs

\begin{document}

\begin{frontmatter}

\title{Improving the Benjamini-Hochberg Procedure for Discrete Tests}
\runtitle{Improving BH for discrete tests}
\begin{aug}
\author{\fnms{Sebastian} \snm{D\"ohler}
\ead[label=e1]{sebastian.doehler@h-da.de}}
\address{
Darmstadt University of Applied Sciences\\ D-64295 Darmstadt, Germany\\
\printead{e1}\\
\phantom{E-mail: sebastian.doehler@h-da.de\ }}
\and
\author{\fnms{Guillermo} \snm{Durand}\ead[label=e2]{guillermo.durand@upmc.fr}}
\address{
Sorbonne Universit\'es, UPMC\\ 4, Place Jussieu, 75252 Paris cedex 05, France\\
\printead{e2}\\
\phantom{Email: guillermo.durand@upmc.fr\ }}
\and
\author{\fnms{Etienne} \snm{Roquain}
\ead[label=e3]{etienne.roquain@upmc.fr}}
\address{
Sorbonne Universit\'es, UPMC\\ 4, Place Jussieu, 75252 Paris cedex 05, France\\
\printead{e3}\\
\phantom{E-mail: etienne.roquain@upmc.fr }}

\end{aug}
\begin{abstract}
	To find interesting items in genome-wide association studies or next generation sequencing data, a crucial point is to design powerful false discovery rate (FDR) controlling procedures that suitably combine discrete tests (typically binomial or Fisher tests). 
	In particular, recent research has been striving for appropriate modifications of the classical  Benjamini-Hochberg (BH) step-up procedure that accommodate discreteness.
	However, despite an important number of attempts, 
	 these procedures did not come with theoretical guarantees.  	  
	The present paper contributes to fill the gap: 
	it presents new modifications of the BH procedure that incorporate the discrete structure of the data and provably control the FDR for any fixed number of null hypotheses (under independence).
	Markedly, our FDR controlling methodology  allows to incorporate simultaneously  the discreteness and 
	the quantity of signal of the data (corresponding therefore to a so-called $\pi_0$-adaptive procedure).  The power advantage of the new methods is demonstrated in a numerical experiment and for some appropriate real data sets.
	\end{abstract}


\end{frontmatter}


Multiple testing procedures are now routinely used to find significant items in massive and complex data. An important focus has been given to methods controlling the false discovery rate (FDR) because this scalable type I error rate ``survives" to high dimension. 
Since the original procedure of 
\cite{BenjaminiHochberg95}, much effort has been undertaken  to design FDR controlling procedures that adapt to various underlying structures of the data, such as the quantity of signal, the signal strength and the dependencies, among others. 

In this work, we deal with adaptation to discrete data. This type of data arises in many relevant applications, in particular when data  are represented by counts. Examples can be found in  clinical studies (see e.g. \cite{WestWolf1997}), genome-wide association studies (GWAS)  (see e.g. \cite{Dickhaus2012}) and next generation sequencing data (NGS)   (see e.g. \cite{DoergeChen2015}). It is well known (see e.g. \cite{WestWolf1997}) that using discrete test statistics can generate a severe power loss, already at the stage of the single tests. A consequence is that using ``blindly" the BH procedure with discrete $p$-values will control the FDR in a  too conservative manner. Therefore, more powerful procedures that avoid this conservatism are much sought after 
 in applications, see for instance \cite{Karp2016}, \cite{vandenBroek2015} and \cite{Dickhaus2012}. 

\begin{figure}[htbp]
	\centering
	\includegraphics[width=0.9\textwidth]{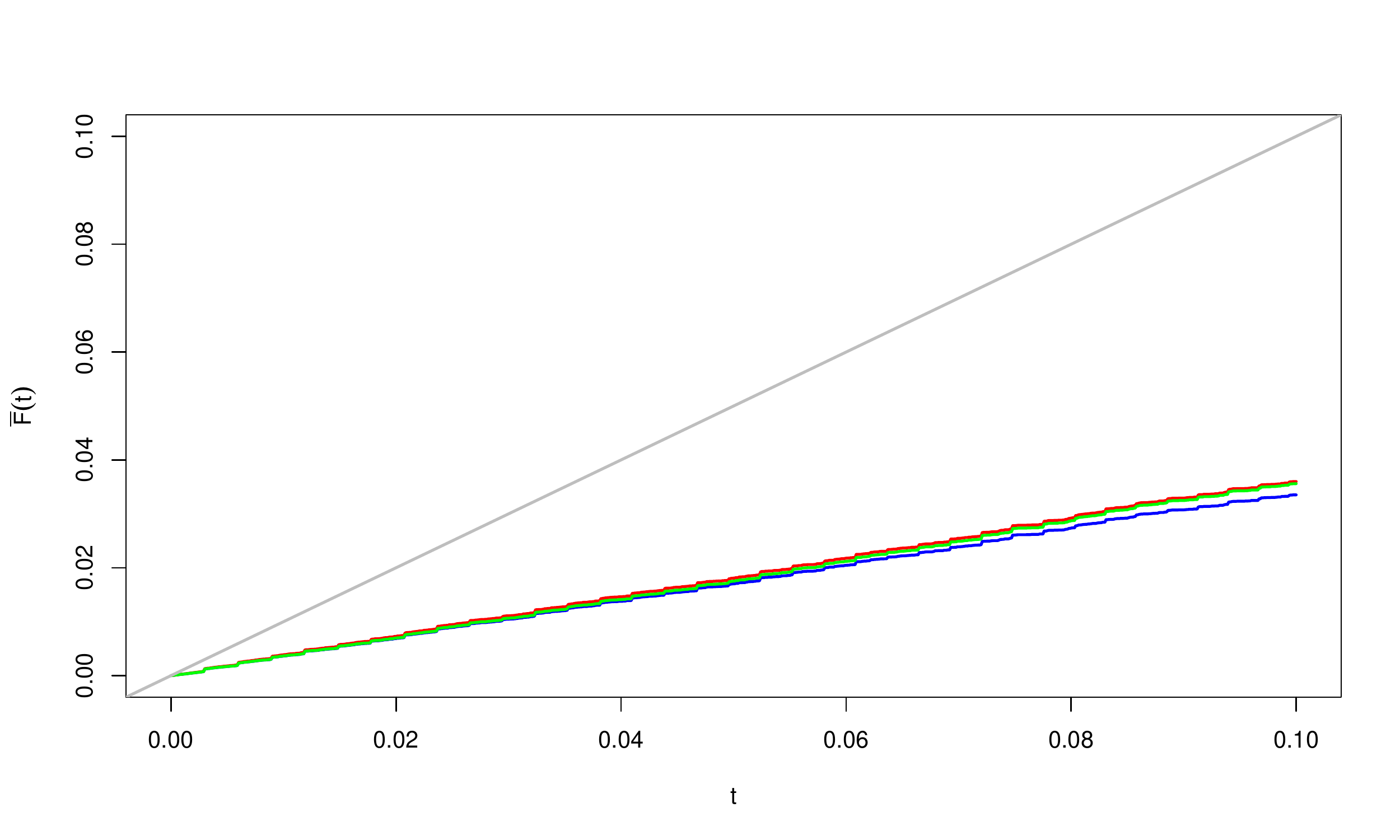}
	\caption{Plots of variants of $\Fbar$ for the pharmacovigilance data. The  grey line corresponds to the uniform case, the discrete variants are represented by blue (for $\Fbar$), green (for $\FbarSD$) and  red (for $\FbarSU$) lines.}
	\label{fig:FbarPlotsHellerForPaper}
\end{figure}

In the literature, building multiple testing procedures that take into account the discreteness of the test statistics has a long history that can be traced back to Tukey and \cite{Mantel1980}.
Some null hypotheses can be {\it a priori} excluded from the study because the corresponding tests are unable to produce sufficiently small $p$-values. This results in a multiplicity reduction that should increase the power. 
While this idea has been exploited in \cite{Tar1990} and in a more general manner in  \cite{WestWolf1997} for family-wise error rate,  
an attempt has been made for FDR later in \cite{Gilbert05}. More recently,
\cite{Heyse2011} has proposed a more powerful solution,  relying on the following averaged cumulative distribution function (c.d.f.):
\begin{equation}\label{equ-Fbar}
	\overline{F}(t) = \frac{1}{m}\sum_{i=1}^m F_i(t), \:\:t\in[0,1],
\end{equation}
where each $F_i$ corresponds to the c.d.f. of the $i$-th test $p$-value under the null hypothesis.
To illustrate the potential benefit of using $\overline{F}$, Figure~\ref{fig:FbarPlotsHellerForPaper} displays this function for the pharmacovigilance data from \citet{Heller2012} (see Section~\ref{sec:EmpiricalData}  for more details). 
The critical values of the Heyse procedure can be obtained by  inverting $\Fbar$ at the values $\alpha k/m, 1\leq k \leq m$. Thus, the smaller the $\Fbar$-values, the larger the critical values.
Here, Heyse critical values improve the BH critical values roughly by a factor 3, thereby yielding a potentially strong rejection enhancement.
Furthermore, since the functions $F_i$'s are known, so is $\Fbar$. Hence, the  user has a good  prior idea of the improvements reachable by this discrete approach.
Unfortunately, the Heyse procedure does not rigorously control the FDR in general; counter-examples are provided in 
\cite{Heller2012} and \cite{Doehler2016}.

Meanwhile, 
different solutions
have been explored by modifying directly the $p$-values, either by randomization (see \cite{Habiger2015} and references therein), or by shrinking them to build so-called {\it midP-values} (see \cite{Heller2012} and references therein). 
Other approaches incorporate discreteness to obtain less conservative FDR estimates, see, e.g., \cite{PoundsCheng06},  or by combining grouping and weighting approaches, see \cite{DoergeChen2015}.


Overall, although many new procedures have been proposed in the literature, only few of them have been proved to achieve a  rigorous FDR control under standard conditions, especially in the finite sample case. 
To the best of our knowledge, we can only refer to the discretised version of the procedure of \cite{BL1999} introduced by \cite{Heller2012} and to the asymptotic work of \cite{Fer2007}.  \cite{DoergeChen2015} sum up the status quo by noting that '{\it \ldots how to derive better FDR procedures in the discrete paradigm remains an urgent but still unresolved problem.}' This paper offers a solution 
by  presenting new procedures that achieve both theoretical validity and good practical performance. \\


The paper is organised as follows: after having precisely defined the setting in Section~\ref{sec:preli}, we introduce in Section~\ref{sec:Methods} new procedures relying on the following modifications of the $\overline{F}$ function
\begin{align*}
	\FbarSU (t)= \frac{1}{m}\sum_{i=1}^m \frac{F_i\left(t\right)}{1-F_i\left(\tau_m\right)};\:\:\:
	\FbarSD (t)=\frac{1}{m}\sum_{i=1}^m \frac{F_i\left(t\right)}{1-F_i\left(t\right)},\:\: t\in[0,1],
\end{align*}
(with the convention $1/0=+\infty$), where an appropriate choice of $\tau_m$ is made.
To feel how light these modifications are, Figure~\ref{fig:FbarPlotsHellerForPaper} displays these functions and shows they are very close to the original $\overline{F}$ for small values of $t$. 
In addition, we also introduce more powerful ``adaptive" versions, meaning that the derived critical values are designed in a way that ``implicitly estimate" the overall proportion of true null hypotheses and thus may outperform the original Heyse procedure. 
Next, in Section~\ref{sec:TheoreticalResults}, we establish rigorous FDR control of the corresponding non-adaptive and adaptive procedures under standard conditions. Our proofs rely on new bounds on FDR that generalise some prominent results of the multiple testing literature. These bounds are the main mathematical contributions of the paper and are interesting in their own right, beyond the discrete setting. 
Also, to explore in detail the improvement of our procedures, 
we analyse both real and simulated data
in Sections~\ref{sec:EmpiricalData} and~\ref{sec:simu}. 
Finally, while the proofs are given in appendix (together with some additional procedures), complementary results are provided in Appendix~\ref{sec:supp}.

%

\section{Preliminaries}\label{sec:preli}

\subsection{General model}\label{sec:model}

Let us observe a random variable $X$, defined on a probabilistic space and valued in an observation space $(\mathcal{X},\mathfrak{X})$. We consider a set $\mathcal{P}$ of possible distributions for the distribution of $X$ and we denote the true one by $P$. 
We assume that $m$ null hypotheses $H_{0,i}$, $1\leq i \leq m$, are available for $P$ and we denote the corresponding set of true null hypotheses by $\cH_0(P)=\{1\leq i \leq m\::\: \mbox{$H_{0,i}$ is satisfied by $P$} \}$. We also denote by $m_0(P)=|\cH_0(P)|$ the number of true nulls.

We assume that the user has at hand a set of $p$-values to test each null, that is, a set of random variables $\{p_i(X), 1\leq i \leq m\}$, valued in $[0,1]$. 
Throughout the paper, we also make the important (but classical) assumption that the  $p$-values $p_i(X)$, $1\leq i \leq m$, are mutually independent.

Now, we denote $\mathcal{F}=\{F_{i}, 1\leq i \leq m\}$, where for each $i\in\{1,\dots,m\}$, the function 
$$
F_{i}(t) = \sup_{P\in \mathcal{P}\::\: i\in  \cH_0(P)} \P_{X\sim P}(p_{i}(X)\leq t ), \:\: t\in[0,1], \:\:1\leq i\leq m
$$
is assumed to be {\it known}. 
Note that we necessarily have $F_i(\cdot)$ non decreasing, $F_{i}(t) \in[0,1]$, $F_i(1)=1$ and we add the technical condition $F_i(0)=0$.
Loosely, each $F_i$ corresponds to the cumulative distribution of $p_i$ under the null.
%
Above, we have taken the supremum to cover the case where the null hypothesis is composite: in that situation, each $F_i$ is adjusted according to the least favorable configuration within the null $H_{0,i}$.

Here are some conditions on $\mathcal{F}$ that will be useful to compare some of the studied procedures (these conditions are {\it not} assumed in our results unless explicitly mentioned):
\begin{align}
&F_{i}(t)\leq t , \:\: t\in[0,1],\:\:1\leq i \leq m, \label{equ-sousunif}\\
&F_{i}(t)= t , \:\: t\in[0,1],\:\:1\leq i \leq m \label{equ-unif}.
\end{align}
Condition \eqref{equ-sousunif} ensures that the $p$-values  have marginals stochastically lower-bounded by a uniform variable under the null, called a {\it super-uniform} distribution in the sequel. This is the classical setting which is used in most of the work dealing with FDR controlling theory, see, e.g., \cite{BenjaminiHochberg95}. 
Condition \eqref{equ-unif} is more restrictive: if each null hypothesis is a singleton, it is equivalent to the $p$-values having uniform marginals under the null.

\subsection{Discrete and continuous modelling}\label{sec:discretecont}
In order to describe the overall support of $p$-value distributions we assume one of the two following situations to be at hand throughout the paper (except in Section~\ref{sec:TheoreticalResults} which is written in a more general manner):
\begin{itemize}
\item Continuous case: for all $i \in\{1,\dots,m\}$, $F_i$ is continuous.  In that case, we let $\mathcal{A}_i=[0,1]$, $1\leq i \leq m$ and $\mathcal{A}=\cup_{i=1}^m \mathcal{A}_i=[0,1]$ is the overall $p$-value support.
\item Discrete case:  each $p$-value $p_i$  (both under the null and alternative) takes values in some finite set $\mathcal{A}_i =\{a_{i,k}, 0\leq k \leq K_i\}$, where $K_i\geq 0$ and $(a_{i,k})_{0\leq k \leq K_i} \in [0,1]^{K_i+1}$ is an increasing sequence (with $a_{i,0}=0$, $a_{i,K_i}=1$). We denote $\mathcal{A}=\cup_{i=1}^m \mathcal{A}_i$ the overall $p$-value support.
\end{itemize}

The continuous setting 
is typically valid in
 situations where the $p$-values are calibrated from test statistics having a continuous distribution under the null. In this situation, \eqref{equ-unif} is often satisfied.
The discrete setting typically arises in situations where the $p$-values are calibrated from test statistics having a finitely supported distribution under the null. In this situation, 
we generally have that \eqref{equ-unif} holds true only on the support of $F_i$, that is,
%
\begin{align}
&F_{i}(t)= t , \:\: t\in\mathcal{A}_i,\:\:1\leq i \leq m. \label{equ-unifgeneral}
\end{align}
In the discrete framework, let us underline that while  \eqref{equ-unifgeneral} will typically hold, the equality $F_i(t)= t$, $t\in\mathcal{A}$ will fail in general because $\mathcal{A}$ contains points of $\mathcal{A}_j$ for $j\neq i$. Then $\overline{F}(t)$ defined by \eqref{equ-Fbar} will be smaller than $t$ in general (see Figure~\ref{fig:FbarPlotsHellerForPaper}), which is exactly the property that we want to exploit in this paper.

To illustrate the above framework, we provide below two simple examples (for more advanced examples, see for instance \cite{DoergeChen2015}). 

\begin{example}[Gaussian testing]
 Observe $X=(X_i)_{1\leq i \leq m}$ with independent coordinates and marginals $X_i\sim \mathcal{N}(\mu_i,1)$, where $\mu_i\in\mathbb{R}$ is the parameter of interest, $1\leq i \leq m$. In that situation, a possible hypothesis testing problem is to consider the nulls $H_{0,i}:$ ``$\mu_i\leq 0$" against $H_{1,i}:$ ``$\mu_i> 0$". Then $p_i(X)=1-\Phi(X_i)$, $1\leq i \leq m$, is a family of $p$-values satisfying \eqref{equ-unif} (where $\Phi$ denotes the c.d.f. of a standard Gaussian variable). 
\end{example}

\begin{example}[Binomial testing]
 Observe $X=(X_i)_{1\leq i \leq m}$ with independent coordinates and marginals $X_i\sim \mathcal{B}(n_i,\theta_i)$, where $n_i\geq 1$ is known and $\theta_i \in (0,1)$ is the parameter of interest, $1\leq i \leq m$. In that situation, a possible hypothesis testing problem is to consider the nulls $H_{0,i}:$ ``$\theta_i\leq 1/2$" against $H_{1,i}:$ ``$\theta_i> 1/2$". Then $p_i(X)=T_i(X_i)$, $1\leq i \leq m$, define a family of $p$-values where $T_i(x)=2^{-n_i} \sum_{j=0}^x {n_i \choose j}$
 is the upper-tail distribution function of a binomial distribution of parameters $(n_i,1/2)$.
 The support of the $p$-values under the null and alternative is covered by letting $K_i=n_i+1$ and $a_{i,k}=2^{-n_i} \sum_{j=0}^{k-1} {n_i \choose j}$, $1\leq k \leq K_i$. 
 We merely check in that case that \eqref{equ-unif} is violated while \eqref{equ-sousunif} and \eqref{equ-unifgeneral} hold. 
\end{example}




\subsection{Step-wise procedures} \label{ssec:StepWiseProcedures}

First define a critical value sequence as any nondecreasing sequence ${\tau}= (\tau_k)_{1\leq k \leq m} \in [0,1]^m$ (with $\tau_0=0$ by convention).

The \textit{step-up procedure} of critical value sequence $\tau$, denoted by ${\SU}(\tau)$,  rejects the $i$-th hypothesis if $p_i\leq \tau_{\hat{k}}$, with 
$\wh{k}= \max\{k\in\{0,1,...,m\}\::\: p_{(k)}\leq {\tau}_k\},$
 where $p_{(1)}\leq p_{(2)} \leq ...\leq p_{(m)}$ denote the ordered $p$-values (with the convention $p_{(0)}=0$). 
 
 The \textit{step-down procedure} of critical value sequence $\tau$, denoted by ${\SD}(\tau)$,  rejects the $i$-th hypothesis if $p_i\leq \tau_{\tilde{k}}$, with 
$\wt{k}= \max\{k\in\{0,1,...,m\}\::\: \forall k'\leq k, \:p_{(k')}\leq {\tau}_{k'}\}.$
 It is straightforward to check that, for the same set of critical values, the step-up version always rejects  more hypotheses than the step-down version. More comments and illustrations on step-wise procedures can be found in \cite{BDRV2014} and \cite{Dic2014}, among others. 

\subsection{False discovery rate}

We measure the quantity of false positives of a step-up (resp. step-down) procedure by using the false discovery rate (FDR), introduced and popularised by \cite{BenjaminiHochberg95}, which is defined as the averaged proportion of errors among the rejected hypotheses. More formally,
for some procedure $R$ rejecting the $i$-th hypothesis if $p_i\leq \hat{t}(X)$ (for some  threshold $\hat{t}(X)$), we let

\begin{equation}\label{eq:FDR}
\FDR(R,P) = \E_{X \sim P}\left[\frac{\sum_{i\in\cH_0(P)} \ind{p_i\leq \hat{t}(X)} }{1\vee \sum_{i=1}^m \ind{p_i\leq \hat{t}(X)} }\right],\:\: P\in\mathcal{P}.
\end{equation}

 The main contribution of this work is to propose procedures that control the FDR at a prescribed level $\alpha$ and that 
incorporate the knowledge of the $F_i$'s in a way that increases the number of discoveries.

\section{Procedures} \label{sec:Methods}
In this section we briefly review some existing methods for FDR control and introduce our new  procedures.
\subsection{Existing methods}
We use the following methods as starting points for constructing new procedures.
\begin{itemize}
\item [-] [BH]:  the seminal procedure proposed in \cite{BenjaminiHochberg95}, corresponding to the step-up procedure ${\SU}(\tau)$, with critical values $\tau_k=\alpha k/m$, $1\leq k \leq m$;
\item [-] [BR-$\lambda$]: an adaptive version of the BH procedure that was proposed in \cite{BR2009}, 
corresponding to the step-up procedure ${\SU}(\tau)$, with critical values 
 \begin{equation}\label{equ-Blacrit}
\tau_k =\left(  (1-\lambda)\frac{ \alpha k}{m-k +1} \right)\wedge \lambda,\:\: 1\leq k \leq m;
\end{equation}
\item [-] [GBS]: an adaptive version of the BH procedure that has been proposed in 
 \cite{Gavrilov2009}, corresponding to the step-down procedure ${\SD}(\tau)$, with critical values 
 \begin{equation}\label{equ-Gavcrit}
\tau_k = \frac{\alpha k}{m-(1-\alpha)k +1},\:\: 1\leq k \leq m;
\end{equation}
\item [-] [Heyse]: 
 the step-up procedure ${\SU}(\tau)$ using critical values given by
\begin{equation}\label{criticalHeyse}
\tau_k = \max\{ t\in \mathcal{A}\::\: \overline{F}(t) \leq \alpha k/m\},\:\: 1\leq k \leq m;
\end{equation}
where $\overline{F}$ is defined by \eqref{equ-Fbar}.
This procedure was proposed in \cite{Heyse2011}.
\end{itemize}

 The rationale behind the critical values of [BR-$\lambda$] and [GBS] is that they are intended to mimic the oracle critical values $\tau_k=\alpha k/m_0(P)$, $1\leq k \leq m$, which are less conservative than those of [BH] when $m_0(P)/m$ is not close to $1$, see, e.g., \cite{BKY2006,BR2009} for more details on adaptive procedures.

Let us now comment on [Heyse]. First, in the continuous setting where \eqref{equ-sousunif} holds, $\overline{F}(t)\leq t$, $t\in[0,1]$, and thus the critical values given by \eqref{criticalHeyse} satisfy $\tau_k \geq \alpha k/m$, $1\leq k \leq m$, which means that [Heyse] rejects at least as many hypotheses as [BH].
When \eqref{equ-unif} additionally holds, we have  $\overline{F}(t)=t$, $t\in[0,1]$, and the two critical value sequences are the same.
Second, in the discrete setting where \eqref{equ-sousunif} holds, $\mathcal{A}$ is finite and $\tau_k$ is not necessarily greater than $ \alpha k/m$ anymore. However, [Heyse] is also less conservative (or equal) than [BH] in the latter case, as stated in the following result (proved in Appendix~\ref{sec:supp} for completeness).

\begin{lemma}\label{lem:HeyseBH}
Consider the model of Section~\ref{sec:model} assuming \eqref{equ-sousunif}, both in the continuous and discrete setting described in Section~\ref{sec:discretecont}. Then the set of nulls rejected by [Heyse] is  larger than the one of [BH] (almost surely). Furthermore, under \eqref{equ-unifgeneral}, these two rejection sets are equal (almost surely) if $F_i=F_j$ for all $i\neq j$.
\end{lemma}

The equality case of Lemma~\ref{lem:HeyseBH} was provided in Proposition~2.3 of \cite{Heller2012}, who presented it as a limitation of Heyse procedure in the discrete case. However, we argue that the condition $F_i=F_j$ for all $i\neq j$ 
is a somehow extreme configuration which is rarely met in practice (in the discrete case). More typically, the $F_i$'s have an heterogeneous structure implying that $\overline{F}(t)$ is smaller than $t$ (see Figure~\ref{fig:FbarPlotsHellerForPaper}). This entails that [Heyse] can substantially improve [BH] (see Figure~\ref{fig:CritConstants}).

While [Heyse] incorporates the knowledge of the $F_i$'s in a natural way (see also Remark~\ref{Rem:bayes} below), it
is not correctly calibrated for a rigorous FDR control: as shown in \cite{Heller2012,Doehler2016}, it fails to control the FDR in general. We propose suitable modifications of [Heyse] in the next sections.

\begin{remark}[Empirical Bayes point of view on the Heyse procedure]\label{Rem:bayes}
We claim that [Heyse] corresponds to a suitable empirical Bayes procedure. 
To see this,  consider the ``binomial example" of Section~\ref{sec:discretecont},   but assume now that the counts $n_1,\dots,n_m$ are observed from a sample $N_1,\dots, N_m$ i.i.d. of a priori distribution $\nu$.  
Unconditionally, the $p$-values $p_i$,  $i\in\cH_0$, are thus i.i.d. with c.d.f.
$
\bar{F}_0=\sum_{n\geq 0} \nu(\{n\}) F_{0,n},
$
where $F_{0,n}$ is the c.d.f. jumping at each $x_{k,n}=2^{-n} \sum_{j=0}^{k-1} {n \choose j}$ with $F_{0,n}(x_{k,n})=x_{k,n}$, $1\leq k \leq n+1$.  
This suggests to normalise the $p$-values $p_i$ as $\bar{F}_0(p_i)$ which leads to the step-up procedure with critical values $\tau_k =\max\{t \::\: \bar{F}_0(t)\leq \alpha k/m\}$. 
Following an empirical Bayes approach, the prior $\nu$ can be estimated by $\hat{\nu}(\{n\}) = m^{-1}\sum_{i=1}^m \mathbf{1}_{\{N_i=n\}}$, which gives rise to the estimator of $\bar{F}_0$ defined  by $
\hat{\bar{F}}_0 = \sum_{n\geq 0} \hat{\nu}(\{n\}) F_{0,n} = m^{-1}\sum_{i=1}^m F_{0,N_i},
 $ 
which is equal to $\overline{F}$ given by \eqref{equ-Fbar}. Hence,
 the corresponding (empirical Bayes) step-up procedure reduces to [Heyse]. 
\end{remark}

\subsection{Two new methods} \label{subsec:NewMethods}

We now present two procedures that aim at correcting [Heyse] :
\begin{itemize}
\item [-] [DBH-SU]: the step-up procedure  ${\SU}(\tau)$ using the critical values defined  in the following way: 
\begin{align}
\tau_m &=\max\left\{ t\in \mathcal{A}\::\: \frac{1}{m}\sum_{i=1}^m \frac{F_i\left(t\right)}{1-F_i\left(t\right)} \leq \alpha \right\}\label{taumDBHSU}\\
\tau_k &= \max\left\{ t\in \mathcal{A}\::\: t\leq \tau_m,\:\frac{1}{m}\sum_{i=1}^m \frac{F_i\left(t\right)}{1-F_i\left(\tau_m\right)} \leq \alpha k/m \right\},\:
1\leq k\leq m-1.\label{critvaluesDBHSU}
\end{align}
\item [-] [DBH-SD]: the step-down procedure  ${\SD}(\tau)$ using the critical values defined  in the following way : 
\begin{align}\label{critvaluesDBHSD}
\tau_k = \max\left\{ t\in \mathcal{A}\::\:  \frac{1}{m}\sum_{i=1}^m \frac{F_i\left(t\right)}{1-F_i\left(t\right)} \leq \alpha k/m \right\},\:1\leq k\leq m.
\end{align}
\end{itemize}

[DBH-SU] can be seen as a correction of [Heyse]: the correction term in the critical values \eqref{critvaluesDBHSU} lies in the additional denominator $1-F_i\left(\tau_m\right)$. A consequence is that [DBH-SU] can be more conservative than [BH]. However, the magnitude of this phenomenon is always small, as the next lemma shows  (proved in Appendix~\ref{sec:supp}  for completeness).

\begin{lemma}\label{lem:DBH-SU-BH}
Under the conditions of Lemma~\ref{lem:HeyseBH},  the set of nulls rejected by [DBH-SU]  contains 
the one of [BH] taken at level $\alpha/(1+\alpha)$ (almost surely). 
\end{lemma}



For [DBH-SD], the following result can be established.

\begin{lemma}\label{lem:DBH-SD-BH}
Under the conditions of Lemma~\ref{lem:HeyseBH},  the set of nulls rejected by [DBH-SD] contains  the one of the step-down procedure with critical values $(\alpha k/m)/(1+\alpha k/m)$, $1\leq k\leq m$ (almost surely). 
\end{lemma}

From \eqref{critvaluesDBHSU} and \eqref{critvaluesDBHSD} it is clear that the critical values of [DBH-SD] are always at least as large as those for [DBH-SU]. However, since the step-up direction is more powerful than the step-down direction (see Section \ref{ssec:StepWiseProcedures}) neither of the two generally dominates the other one.

\begin{remark}
We may ask whether we can construct a uniform improvement of [BH] that incorporate the $F_i$'s. There is indeed such a procedure (see procedure [RBH] in Appendix~\ref{sec:RBH} for more details). However, the improvement brought by the $F_i$'s information is less substantial than for [DBH-SU], so we have chosen to omit [RBH] from the main stream of the paper. 

\end{remark}

\subsection{Adaptive versions}

In this section, we define adaptive versions of [DBH-SU] and [DBH-SD] in the following way:
\begin{itemize}
\item [-] [A-DBH-SU]: the step-up procedure  ${\SU}(\tau)$ using 
the critical values defined  in the following way: $\tau_m$ as in \eqref{taumDBHSU} and for $1\leq k\leq m-1$,
\begin{align}
\tau_k &= \max\left\{ t\in \mathcal{A}\::\: t\leq \tau_m,\:  \left(  \frac{F\left(t\right)}{1-F\left(\tau_m\right)}\right)_{(1)}+\dots + \left(  \frac{F\left(t\right)}{1-F\left(\tau_m\right)}\right)_{(m-k+1)}\leq \alpha k\right\},\:
\label{critvaluesADBHSU}
\end{align}
 where each $\left( \frac{F\left(t\right)}{1-F\left(\tau_m\right)}\right)_{(j)}$ denotes the $j$-th largest elements of the set $\left\{ \frac{F_i\left(t\right)}{1-F_i\left(\tau_m\right)} , 1\leq i \leq m \right\}$.

\item [-] [A-DBH-SD]: the step-down procedure  ${\SD}(\tau)$ using the critical values defined  in the following way : 
\begin{align}\label{critvaluesADBHSD}
\tau_k = \max\left\{ t\in \mathcal{A}\::\:   \left(  \frac{F\left(t\right)}{1-F\left(t\right)}\right)_{(1)}+\dots + \left(  \frac{F\left(t\right)}{1-F\left(t\right)}\right)_{(m-k+1)}\leq \alpha k\right\},\:1\leq k\leq m,
\end{align}
where each $\left( \frac{F\left(t\right)}{1-F\left(t\right)}\right)_{(j)}$ denotes the $j$-th largest elements of the set $\left\{ \frac{F_i\left(t\right)}{1-F_i\left(t\right)} , 1\leq i \leq m \right\}$.
\end{itemize}

Note that the critical values of [A-DBH-SU]  and [A-DBH-SD] are clearly larger than or equal to those of their non-adaptive counterparts [DBH-SU]  and [DBH-SD], respectively. This means that the adaptive versions are always less conservative. 

The following result establishes a connection of the adaptive procedures to the [BR-$\lambda$] and [GBS] procedures (proved in Appendix~\ref{sec:supp} for completeness). 

\begin{lemma}\label{lem:adapt}
Under the conditions of Lemma~\ref{lem:HeyseBH}, the following holds:
\begin{itemize}
\item[(i)] the set of nulls rejected by [A-DBH-SU]  contains the one of  [BR-$\lambda$] (almost surely), where $\lambda$ is taken equal to \eqref{taumDBHSU};
\item[(ii)]  the set of nulls rejected by [A-DBH-SD]  contains the one of [GBS] (almost surely); 
\end{itemize}
\end{lemma}

The above lemma ensures that the user can incorporate the knowledge of the $F_i$'s  in adaptive procedures with a ``no loss" guarantee with respect to [BR] and [GBS]. This is a somehow striking fact, coming loosely from a ``fortunate marriage" between the proof technics of discreteness theory and adaptation theory. 

\begin{remark}
We may ask whether we can build a procedure that is a uniform improvement of [BR-$\lambda$], for any fixed value of $\lambda\in(0,1)$. We propose a solution in Appendix~\ref{sec:DBR}, called [DBR-$\lambda$]. It does not improve uniformly [DBH-SU], but is an interesting variant of [A-DBH-SU].
\end{remark}

\section{New  FDR bounds} \label{sec:TheoreticalResults}


In this section, we present new FDR bounds which are the main mathematical contributions of this paper and that are of independent interest.
They generalise some classical bounds from super-uniform null distributions to arbitrary heterogeneous (not necessarily discrete) null distributions, and immediately yield FDR control of our new procedures.

\subsection{Results}

First, remember that the model of Section~\ref{sec:model} basically only assume independence between the $p$-values (and not super-uniformity of the null distribution). The following result holds.

\begin{theorem}\label{th:step-updown}
In the model of Section~\ref{sec:model}, for any critical values $\tau_{k}$, $1\leq k \leq m$ and for all  $P\in\mathcal{P}$, we have
\begin{align}
\FDR(\SU(\tau),P) &\leq\min\left(\sum_{i=1}^m  \max_{1\leq k \leq m} \frac{F_i(\tau_k)}{k},
\max_{1\leq k \leq m}\max_{ \substack{A \subset \{1,\dots,m\}\\ |A|=m-k+1} }\left( \frac{1}{k}  \sum_{i\in A} \frac{F_i\left({\tau}_{k}\right)}{1-F_i\left({\tau}_{m}\right)} \right)
\right);\label{SUFDRbound}\\
\FDR(\SD(\tau),P) &\leq \min\left(\sum_{i=1}^m  \max_{1\leq k \leq m} \frac{F_i(\tau_k)}{k},
 \max_{1\leq k \leq m}\max_{ \substack{A \subset \{1,\dots,m\}\\ |A|=m-k+1} }\left( \frac{1}{k}  \sum_{i\in A} \frac{F_i\left({\tau}_{k}\right)}{1-F_i\left({\tau}_{k}\right)} \right)  \right).\label{SDFDRbound}
\end{align}
\end{theorem}

The proof of Theorem~\ref{th:step-updown} is deferred to Appendix~\ref{sec:proofs}. 
It combines   several techniques: the first tool is an expression of the FDR  introduced by \cite{Fer2007} (step-up case) and \cite{RV2011} (step-down case). 
A second idea comes from the work  \cite{BR2009} (step-up case) and \cite{Gavrilov2009} (step-down case), which introduced a new term (here,  the denominator $(1-F_i(\cdot))$) to make the adaptive argument works fine. Finally, another inspiration is the study of \cite{RW2009} and \cite{Doehler2016}  that allowed to deal with  heterogeneous FDR thresholding.
Let us underline that  the obtained proof is especially concise, which means that these different techniques fit together perfectly well, which is perhaps surprising at first glance, see Appendix~\ref{sec:proofsd}. 


Next, let us note that taking the maximum over the subset $A$ in  \eqref{SUFDRbound} and \eqref{SDFDRbound} allows us to adapt to the unknown number of true null hypotheses: loosely, if $k-1$ is the number of rejections, $A$ corresponds to the acceptation set (hence of cardinality $m-k+1$), which ``estimates" $\cH_0$ and thus the sums in   \eqref{SUFDRbound} and \eqref{SDFDRbound} are indexed by a set ``close" to the unknown set $\cH_0$. 
Taking the maximum then corresponds to account for the least favorable possible $\cH_0$. 

Finally, let us underline again that the above bounds do not use the super-uniformity of the $F_i$'s which makes them  quite general and flexible tools. As a case in point, consider mid-$p$-values which were introduced by \cite{Lancaster1961} and are sometimes  used for analysing discrete data (see e.g. \cite{Karp2016}). These $p$-values are no longer super-uniform under the null hypotheses, however our theorem can accomodate such distributions in a natural way to still yield valid FDR controlling procedures. In addition, note that our bounds can be useful outside the discrete setting, when the $F_i$'s are continuous but with flat parts, see the (toy) Example~\ref{rem:toyexample} below.


\subsection{Rationale and relation to  previous work}

Let us now give some intuition behind these bounds by showing how it allows to cover previous work in the literature. 

First, assuming the super-uniformity $F_i(t)\leq t$ for all $t$ and $i$, then these bounds entail 
 \begin{align}
 \FDR(\SU(\tau),P) &\leq m \max_{1\leq k \leq m}\{\tau_k /k\};\label{recover1} \\
\FDR(\SU(\tau),P) &\leq 
\max_{1\leq k \leq m} \frac{m-k+1}{1-{\tau}_{m}} \frac{{\tau}_{k}}{k} ;\label{recover2}\\
\FDR(\SD(\tau),P) &\leq \max_{1\leq k \leq m} \frac{m-k+1}{1-{\tau}_{k}} \frac{{\tau}_{k}}{k},\label{recover3}
\end{align}
which immediately recover the fact that [BH], [BR-$\lambda$] (with $\tau_m=\lambda$) and [GBS] all control the FDR at level $\alpha$. To this respect, bounds \eqref{recover1}, \eqref{recover2} and \eqref{recover3} 
 encompass Theorem~\ref{th:step-updown} of \cite{BenjaminiHochberg95}, Theorem~9 of \cite{BR2009} and Theorem~1.1 of \cite{Gavrilov2009}, respectively. 

Second, by removing the adaptative part of the bounds, that is, by replacing $A$ by $\{1,\dots,m\}$, we obtain the simpler but more conservative bounds
\begin{align}
\FDR(\SU(\tau),P) &\leq 
\max_{1\leq k \leq m}\left( \frac{1}{k}  \sum_{i=1}^m \frac{F_i\left({\tau}_{k}\right)}{1-F_i\left({\tau}_{m}\right)} \right) =\max_{1\leq k \leq m}  m\FbarSU (\tau_k)/k;\label{FDRSUboundsnonadapt}\\
\FDR(\SD(\tau),P) &\leq 
 \max_{1\leq k \leq m}\left( \frac{1}{k}  \sum_{i=1}^m \frac{F_i\left({\tau}_{k}\right)}{1-F_i\left({\tau}_{k}\right)} \right) = \max_{1\leq k \leq m} m\FbarSD (\tau_k)/k\label{FDRSDboundsnonadapt},
\end{align}
where $\FbarSU $ and $\FbarSD $ are defined in the introduction, see also Figure~\ref{fig:FbarPlotsHellerForPaper}.
These variants illustrate perhaps more intuitively how the Heyse-type procedures take advantage of the heterogeneous structure: if some of the $F_i$'s are really small, they will not contribute much into $\FbarSU $ (or $\FbarSD$), offering some additional room  for the other $F_j$'s. 

Finally, these bounds immediately imply that our new procedures enjoy the desired FDR controlling property.

\begin{corollary}\label{cor}
In the model of Section~\ref{sec:model},  both in the continuous and discrete setting described in Section~\ref{sec:discretecont}, 
the procedures [DBH-SU]; [DBH-SD]; [A-DBH-SU]; [A-DBH-SD] all control the FDR at level $\alpha$.
 \end{corollary}

%



\begin{example}\label{rem:toyexample}
Assume that 
the hypotheses are structured in $3$ non-overlapping groups $S_1$, $S_2$ and $S_3$, each of cardinality $m/3$ (assumed to be an integer). Assume that $F_i(x)$ is equal to $x$ if $i\in S_1$, $0$ if $i\in S_2$, and $F_i(x)=2x $ ($x\in[0,1/4]$); $1/2$ ($x\in[1/4,3/4]$); $2x-1$ ($x\in[3/4,1]$), if $i\in S_3$. Then the bound \eqref{FDRSDboundsnonadapt} becomes (for $\tau_k\leq 1/4$):
$$
 \max_{1\leq k \leq m} \frac{m \tau_k}{3 k} \left( \frac{1}{1-\tau_k} +\frac{2}{1-2\tau_k}\right), 
$$
 which entails a new step-down FDR controlling procedure by taking $\tau_k$ such that the above expression is equal to $\alpha$. For $\alpha$ small, we get $\tau_k \approx \alpha k /m$ which yields to a procedure close to a step-down version of [BH]. It thus controls the FDR even though the super-uniformity of the $p$-values is violated. In particular, this illustrates that our methodology exceeds the scope of discrete tests.
\end{example}

\section{Empirical data} \label{sec:EmpiricalData}
 To illustrate the performance of FDR-controlling  procedures for discrete data, we analyse two benchmark data sets which have also been used in previous publications. In what follows, our main goal is to compare the performance of the new procedures [DBH-SU], [A-DBH-SU] and [A-DBH-SD]   to the classical [BH] procedure. As a further benchmark we also include [Heyse] in the analysis. All analyses were performed using the R language for statistical computing \citep{RSoft1}.

\subsection{Pharmacovigilance data}
This data set is derived from a database for reporting, investigating and monitoring adverse drug reactions due to the Medicines and Healthcare products Regulatory Agency in the United Kingdom. It contains the number of reported cases of amnesia as well as the total number of adverse events reported for each of the $m=2446$ drugs in the database. For more details we refer to \citet{Heller2012} and to the accompanying R-package 'discreteMTP' \citep{discreteMTP},  which also contains the data. \citet{Heller2012} investigate the association between reports of amnesia and suspected drugs by performing for each drug a Fisher's exact test (one-sided) for testing association between the drug and amnesia while adjusting for multiplicity by using several (discrete) FDR procedures.

\subsection{Next generation sequencing data }
We also revisit the next generation sequencing (NGS) count data analysed by \citet{DoergeChen2015}, to which we also refer for more details. More specifically, we reanalyse the methylation data set for cytosines of Arabidopsis in \citet{Lister2008} which is part of the R-package 'fdrDiscreteNull' \citep{fdrDiscreteNull}. This data set contains the counts for a biological entity under two different biological conditions or treatments. Following \citet{DoergeChen2015}, $m=7421$ genes whose treatment-wise total counts are positive but row-total counts are no greater than 100 are analysed using the exact binomial test, see \citet{DoergeChen2015}.

\subsection{Results }
Table \ref{tab:Empirical data} summarises the number of discoveries for the pharmacovigilance and NGS data when using the respective FDR procedures at level  $\alpha=0.05$. 
\begin{table}
\caption{Number of rejections (discoveries) for the pharmacovigilance and Arabidopsis methylation data.\label{tab:Empirical data}}
\centering
	\begin{tabular}{lccccc}
		\toprule
		\multicolumn{1}{c}{Data set} & [BH]    & [DBH-SU] 
		 & [Heyse] & [A-DBH-SU] 
		 & [A-DBH-SD] 
		  \\
		\midrule
		Pharmacovigilance  & 24    & 27    & 27    & 27    & 27 \\
		Arabidopsis methylation & 2097  & 2358  & 2379  & 2446  & 2453 \\
		\bottomrule
	\end{tabular}%
\end{table}%
Compared to the classical FDR controlling procedures, the new procedures are able to detect three additional candidates linking amnesia and drugs  in the pharmacovigilance data. Note also that for this data, they reject the same number of hypotheses as [Heyse], even though [Heyse] is not correctly calibrated for FDR control. For the Arabidopsis data, the new procedures improve considerably on [BH]. Moreover, there is a clear separation between the adaptive and non-adaptive procedures.

\begin{figure}[htbp]
	\centering
	\makebox{\includegraphics[width=1\linewidth]{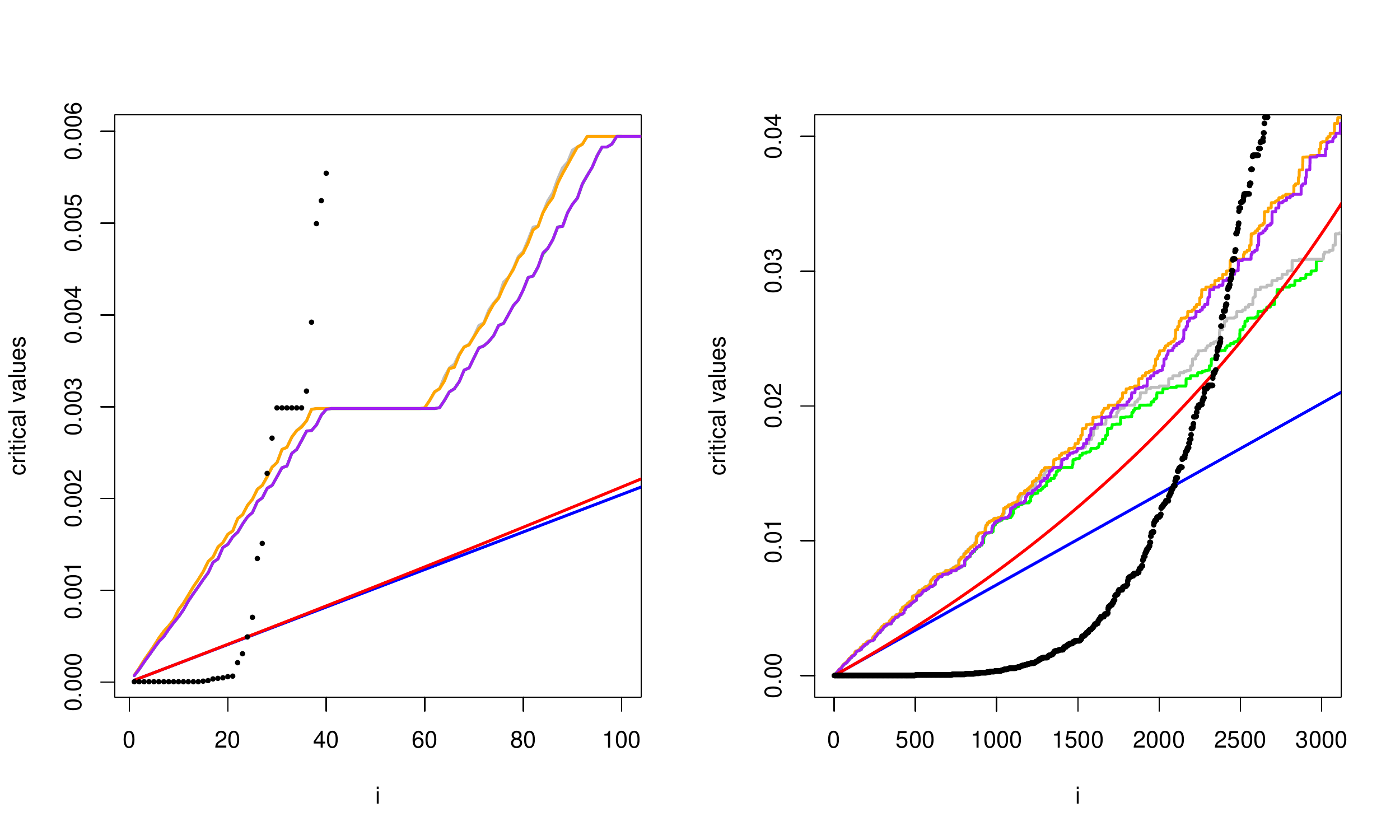}}
	\caption{Critical constants and sorted $p$-values (represented by black dots) for the pharmacovigilance (left panel) and Arabidopsis methylation data (right panel). The [BH], [DBH-SU], [A-DBH-SU], [A-DBH-SD]  and [Heyse] critical constants are represented respectively by blue, green, purple, orange and grey solid lines.	\label{fig:CritConstants}}
\end{figure}

Figure \ref{fig:CritConstants} illustrates graphically the data and  the critical  constants of the involved multiple testing procedures.
In particular, the benefit of taking discreteness into account becomes more apparent: for the pharmacovigilance data, the discrete critical values are considerably (by a factor of $2.5 -3.5$) larger than their respective classical counterparts. This leads to more powerful procedures. For the NGS data, we can observe quite clearly that the [DBH-SU] critical constants are dominated by the [A-DBH-SU] constants, as explained  in Section~\ref{sec:Methods}. This leads to roughly 100 additional rejections. Again, the discrete critical values are  considerably larger than their respective classical counterparts. In \ref{subsec:NewMethods} we mentioned that the correction factor $1- F_i(\tau_m)$, introduced for guaranteeing FDR control of [DBH-SU], may lead to a procedure which is more conservative than [BH]. However, Figure \ref{fig:CritConstants}  shows that -- at least for the data sets considered here -- this risk is by far compensated by the benefit of taking discreteness adequately into account.

\section{Simulation study}\label{sec:simu}
We now investigate the power of the procedures from the previous section in a simulation study similar to those described in \cite{Gilbert05}, \cite{Heller2012} and \cite{Doehler2016}. Again, we focus on comparing  the performance of the new discrete  procedures to [BH].



\subsection{Simulated Scenarios}

We simulate a two-sample problem in which a vector of $m$ independent binary responses (``adverse events") is observed for each subject in two groups, where each group consists of $N=25$ subjects. Then, the goal is to simultaneously test the $m$ null hypotheses $H_{0i}:$ ``$p_{1i}=p_{2i}$", $i=1,\ldots,m$, where $p_{1i}$ and $p_{2i}$ are the success probabilities for the $i$th binary response in group 1 and 2, respectively. We take $m=800,2000$ where $m=m_{1}+m_{2}+m_{3}$ and data are generated so that the response is $Bernoulli(0.01)$ at $m_{1}$ positions for both groups, $Bernoulli(0.10)$ at $m_{2}$ positions for both groups and   $Bernoulli(0.10)$ at $m_{3}$ positions for group 1 and  $Bernoulli(q)$ at $m_{3}$ positions for group 2 where $q=0.15,0.25,0.4$ represents weak,  moderate and strong effects respectively. 
The null hypothesis is true for the  $m_{1}$ and $m_{2}$ positions while the alternative hypothesis is true for the  $m_{3}$ positions. We also take different configurations for the proportion of false null hypotheses, $m_{3}$ is set to be $10\%$,  $30\%$ and  $80\%$ of the value of $m$, which represents  small, intermediate and large proportion of effects (the proportion of true nulls $\pi_{0}$ is $0.9$, $0.7$, $0.2$, respectively). Then, $m_{1}$ is set to be $20\%$,  $50\%$ and  $80\%$ of the number of true nulls (that is, $m-m_{3}$) and $m_{2}$ is taken accordingly as $m-m_{1}-m_{3}$. 

For each of the 54 possible parameter configurations specified by $m,m_{3},m_{1}$ and $q$, $10000$ Monte Carlo trials are performed, that is, $10000$ data sets are generated and for each data set, an unadjusted two-sided $p$-value from Fisher's exact test is computed for each of the $m$ positions, and the multiple testing procedures mentioned above are applied at level $\alpha=0.05$. The power of each procedure was estimated as the fraction of the $m_{3}$ false null hypotheses that were rejected, averaged over the $10000$ simulations. For random number generation the R-function \textit{rbinom} was used. The two-sided $p$-values from Fisher's exact test were computed using the R-function \textit{fisher.test}.

\subsection{Results}

\begin{figure}[htbp]
	\centering
		\makebox{\includegraphics[width=1.0\textwidth]{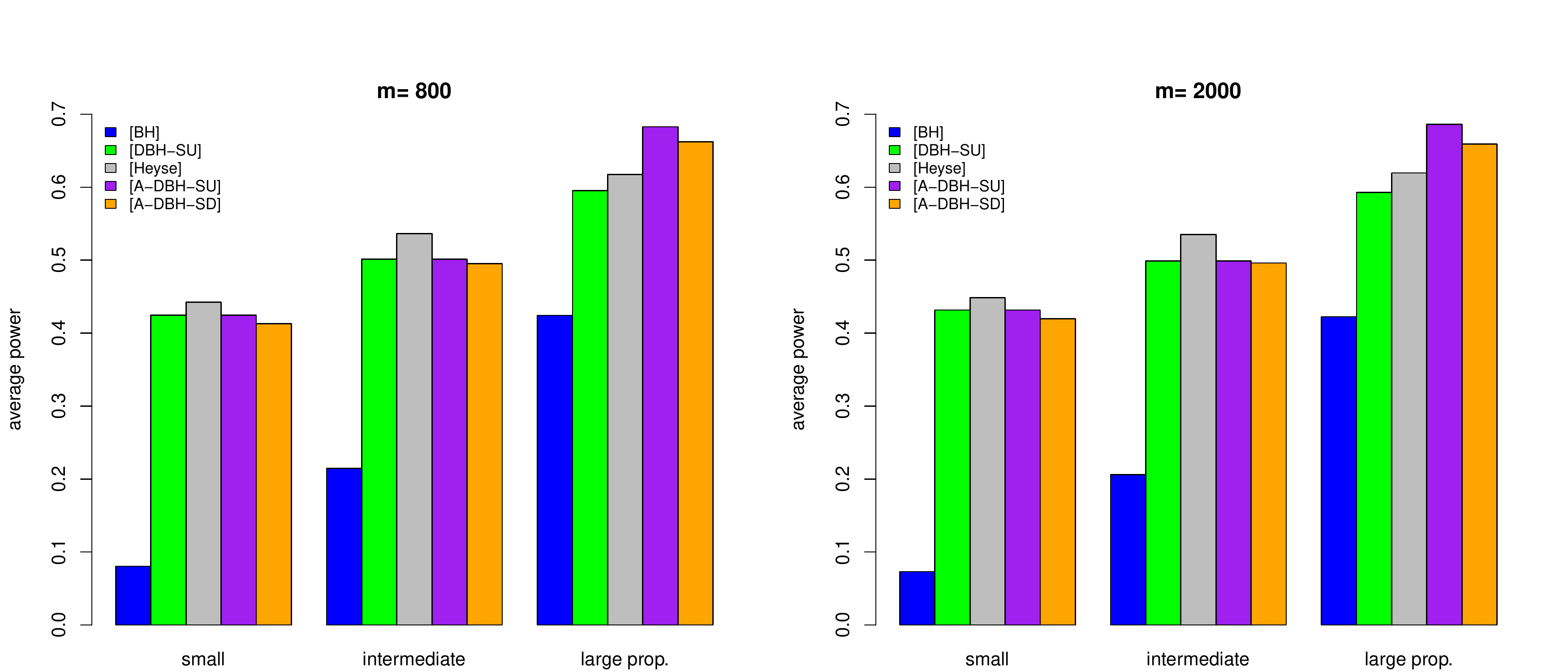}}
	\caption{Average power for the [BH] and discrete procedures in the simulation study. The coloring is the same as in Figure \ref{fig:CritConstants}}
	\label{fig:FDPvsModifiedFDPAverageProportionRejections_25}
\end{figure}

We have computed the (average) power of the five procedures under investigation in all the scenarios (see Table~1 in Appendix~\ref{sec:supp} for the full display).  For weak and moderate effects, i.e. $q=0.15$ and $q=0.25$, none of the procedure possesses relevant power.  For strong effects, the results are summarised 
  in Figure \ref{fig:FDPvsModifiedFDPAverageProportionRejections_25}. (Since 
the power of the discrete procedures is slightly increasing in $m_1$ for fixed $m_3$ and $q$, we present -- in order to avoid over-optimism -- the configuration with smallest $m_1$).


The results are consistent with the findings of the previous section: The new discrete procedures are considerably more powerful than their classical counterparts and perform roughly similarly for small and intermediate proportions of alternatives. When the proportion of alternatives is large, the benefit of using adaptive procedures -- especially the [A-DBH-SU] procedure --  is clearly visible in Figure \ref{fig:FDPvsModifiedFDPAverageProportionRejections_25}. 

\section{Conclusion and discussion}

In this paper, we provided new bounds for the FDR of step-up and step-down procedures that use discrete test statistics. This allowed to define a new class of multiple testing procedures that provably control the FDR while they incorporated  the discreteness of the tests statistics in a convenient way.
We have shown that
our approach can be seen as correcting and improving the approach of \cite{Heyse2011}: while it ensures a theoretical control, it can also make more rejections when the signal amplitude is strong enough.

Let us also mention that 
our procedures are fully usable in practice. 
An R-package is in preparation and it will be presented in a companion paper (which will also deal with computational aspects).

Finally, this paper opens several directions for future research, especially by trying to extend our arguments to other frameworks. For instance, it may be worth to relax the independence requirements. To this respect, we believe that our procedures will inheritate the behavior of BH procedure: while the FDR control is likely to be  maintained  under ``realistic" dependence, formally proving such a result is probably a challenging problem. 


\section*{Acknowledgements}

This work has been supported by the CNRS (PEPS FaSciDo) and the French grant ANR-16-CE40-0019.

\bibliographystyle{chicago}
\bibliography{biblio}

\begin{thebibliography}{}

\bibitem[\protect\citeauthoryear{Benjamini and Hochberg}{Benjamini and
  Hochberg}{1995}]{BenjaminiHochberg95}
Benjamini, Y. and Y.~Hochberg (1995).
\newblock Controlling the false discovery rate: A practical and powerful
  approach to multiple testing.
\newblock {\em Journal of the Royal Statistical Society. Series B\/}~{\em
  57\/}(1), 289--300.

\bibitem[\protect\citeauthoryear{Benjamini, Krieger, and Yekutieli}{Benjamini
  et~al.}{2006}]{BKY2006}
Benjamini, Y., A.~M. Krieger, and D.~Yekutieli (2006).
\newblock Adaptive linear step-up procedures that control the false discovery
  rate.
\newblock {\em Biometrika\/}~{\em 93\/}(3), 491--507.

\bibitem[\protect\citeauthoryear{Benjamini and Liu}{Benjamini and
  Liu}{1999}]{BL1999}
Benjamini, Y. and W.~Liu (1999).
\newblock A step-down multiple hypotheses testing procedure that controls the
  false discovery rate under independence.
\newblock {\em J. Statist. Plann. Inference\/}~{\em 82\/}(1-2), 163--170.

\bibitem[\protect\citeauthoryear{{Blanchard}, {Dickhaus}, {Roquain}, and
  {Villers}}{{Blanchard} et~al.}{2014}]{BDRV2014}
{Blanchard}, G., T.~{Dickhaus}, E.~{Roquain}, and F.~{Villers} (2014).
\newblock {On least favorable configurations for step-up-down tests}.
\newblock {\em Statist. Sinica\/}~{\em 24\/}(1), 1--23.

\bibitem[\protect\citeauthoryear{Blanchard and Roquain}{Blanchard and
  Roquain}{2009}]{BR2009}
Blanchard, G. and E.~Roquain (2009).
\newblock Adaptive false discovery rate control under independence and
  dependence.
\newblock {\em J. Mach. Learn. Res.\/}~{\em 10}, 2837--2871.

\bibitem[\protect\citeauthoryear{Chen and Doerge}{Chen and
  Doerge}{2015a}]{fdrDiscreteNull}
Chen, X. and R.~Doerge (2015a).
\newblock {\em fdrDiscreteNull: False Discovery Rate Procedure Under Discrete
  Null Distributions}.
\newblock R package version 1.0.

\bibitem[\protect\citeauthoryear{Chen and Doerge}{Chen and
  Doerge}{2015b}]{DoergeChen2015}
Chen, X. and R.~Doerge (2015b).
\newblock A weighted fdr procedure under discrete and heterogeneous null
  distributions.
\newblock {\em arXiv:1502.00973\/}.

\bibitem[\protect\citeauthoryear{Dickhaus}{Dickhaus}{2014}]{Dic2014}
Dickhaus, T. (2014).
\newblock {\em Simultaneous statistical inference}.
\newblock Springer, Heidelberg.
\newblock With applications in the life sciences.

\bibitem[\protect\citeauthoryear{Dickhaus, Stra{\ss}burger, Schunk,
  Morcillo-Suarez, Illig, and Navarro}{Dickhaus et~al.}{2012}]{Dickhaus2012}
Dickhaus, T., K.~Stra{\ss}burger, D.~Schunk, C.~Morcillo-Suarez, T.~Illig, and
  A.~Navarro (2012).
\newblock {How to analyze many contingency tables simultaneously in genetic
  association studies.}
\newblock {\em Statistical applications in genetics and molecular
  biology\/}~{\em 11\/}(4).

\bibitem[\protect\citeauthoryear{{D{\"o}hler}}{{D{\"o}hler}}{2016}]{Doehler2016}
{D{\"o}hler}, S. (2016).
\newblock A discrete modification of the {B}enjamini{–-}{Y}ekutieli
  procedure.
\newblock {\em Econometrics and Statistics\/}.

\bibitem[\protect\citeauthoryear{Ferreira}{Ferreira}{2007}]{Fer2007}
Ferreira, J.~A. (2007).
\newblock The {B}enjamini-{H}ochberg method in the case of discrete test
  statistics.
\newblock {\em Int. J. Biostat.\/}~{\em 3}, Art. 11, 18.

\bibitem[\protect\citeauthoryear{Ferreira and Zwinderman}{Ferreira and
  Zwinderman}{2006}]{FZ2006}
Ferreira, J.~A. and A.~H. Zwinderman (2006).
\newblock On the {B}enjamini-{H}ochberg method.
\newblock {\em Ann. Statist.\/}~{\em 34\/}(4), 1827--1849.

\bibitem[\protect\citeauthoryear{Gavrilov, Benjamini, and Sarkar}{Gavrilov
  et~al.}{2009}]{Gavrilov2009}
Gavrilov, Y., Y.~Benjamini, and S.~K. Sarkar (2009).
\newblock An adaptive step-down procedure with proven {FDR} control under
  independence.
\newblock {\em Ann. Statist.\/}~{\em 37\/}(2), 619--629.

\bibitem[\protect\citeauthoryear{Gilbert}{Gilbert}{2005}]{Gilbert05}
Gilbert, P. (2005).
\newblock A modified false discovery rate multiple-comparisons procedure for
  discrete data, applied to human immunodeficiency virus genetics.
\newblock {\em Journal of the Royal Statistical Society. Series C\/}~{\em
  54\/}(1), 143--158.

\bibitem[\protect\citeauthoryear{Habiger}{Habiger}{2015}]{Habiger2015}
Habiger, J.~D. (2015).
\newblock Multiple test functions and adjusted {$p$}-values for test statistics
  with discrete distributions.
\newblock {\em J. Statist. Plann. Inference\/}~{\em 167}, 1--13.

\bibitem[\protect\citeauthoryear{{Heller} and {Gur}}{{Heller} and
  {Gur}}{2011}]{Heller2012}
{Heller}, R. and H.~{Gur} (2011, December).
\newblock {False discovery rate controlling procedures for discrete tests}.
\newblock {\em ArXiv e-prints\/}.

\bibitem[\protect\citeauthoryear{Heller, Gur, and Yaacoby}{Heller
  et~al.}{2012}]{discreteMTP}
Heller, R., H.~Gur, and S.~Yaacoby (2012).
\newblock {\em discreteMTP: Multiple testing procedures for discrete test
  statistics}.
\newblock R package version 0.1-2.

\bibitem[\protect\citeauthoryear{Heyse}{Heyse}{2011}]{Heyse2011}
Heyse, J.~F. (2011).
\newblock A false discovery rate procedure for categorical data.
\newblock In {\em Recent Advances in Bio- statistics: False Discovery Rates,
  Survival Analysis, and Related Topics}, pp.\  43--58.

\bibitem[\protect\citeauthoryear{Karp, Heller, Yaacoby, White, and
  Benjamini}{Karp et~al.}{2016}]{Karp2016}
Karp, N.~A., R.~Heller, S.~Yaacoby, J.~K. White, and Y.~Benjamini (2016).
\newblock Improving the identification of phenotypic abnormalities and sexual
  dimorphism in mice when studying rare event categorical characteristics.
\newblock {\em Genetics\/}.

\bibitem[\protect\citeauthoryear{Lancaster}{Lancaster}{1961}]{Lancaster1961}
Lancaster, H.~O. (1961).
\newblock Significance tests in discrete distributions.
\newblock {\em Journal of the American Statistical Association\/}~{\em
  56\/}(294), 223--234.

\bibitem[\protect\citeauthoryear{Lister, O'Malley, Tonti-Filippini, Gregory,
  Berry, Millar, and Ecker}{Lister et~al.}{2008}]{Lister2008}
Lister, R., R.~C. O'Malley, J.~Tonti-Filippini, B.~D. Gregory, C.~C. Berry,
  A.~H. Millar, and J.~R. Ecker (2008, May).
\newblock {Highly integrated single-base resolution maps of the epigenome in
  Arabidopsis.}
\newblock {\em Cell\/}~{\em 133\/}(3), 523--536.

\bibitem[\protect\citeauthoryear{Mantel}{Mantel}{1980}]{Mantel1980}
Mantel, N. (1980).
\newblock A biometrics invited paper. assessing laboratory evidence for
  neoplastic activity.
\newblock {\em Biometrics\/}~{\em 36\/}(3), 381--399.

\bibitem[\protect\citeauthoryear{Pounds and Cheng}{Pounds and
  Cheng}{2006}]{PoundsCheng06}
Pounds, S. and C.~Cheng (2006).
\newblock {Robust estimation of the false discovery rate}.
\newblock {\em Bioinformatics\/}~{\em 22\/}(16), 1979--1987.

\bibitem[\protect\citeauthoryear{{R Core Team}}{{R Core Team}}{2016}]{RSoft1}
{R Core Team} (2016).
\newblock {\em R: A Language and Environment for Statistical Computing}.
\newblock Vienna, Austria: R Foundation for Statistical Computing.

\bibitem[\protect\citeauthoryear{Roquain and van~de Wiel}{Roquain and van~de
  Wiel}{2009}]{RW2009}
Roquain, E. and M.~van~de Wiel (2009).
\newblock Optimal weighting for false discovery rate control.
\newblock {\em Electron. J. Stat.\/}~{\em 3}, 678--711.

\bibitem[\protect\citeauthoryear{Roquain and Villers}{Roquain and
  Villers}{2011}]{RV2011}
Roquain, E. and F.~Villers (2011).
\newblock Exact calculations for false discovery proportion with application to
  least favorable configurations.
\newblock {\em Ann. Statist.\/}~{\em 39\/}(1), 584--612.

\bibitem[\protect\citeauthoryear{Tarone}{Tarone}{1990}]{Tar1990}
Tarone, R.~E. (1990).
\newblock A modified bonferroni method for discrete data.
\newblock {\em Biometrics\/}~{\em 46\/}(2), 515--522.

\bibitem[\protect\citeauthoryear{van~den Broek, Dijkstra, Krijgsman, Sie, Haan,
  Traets, van~de Wiel, Nagtegaal, Punt, Carvalho, Ylstra, Abeln, Meijer, and
  Fijneman}{van~den Broek et~al.}{2015}]{vandenBroek2015}
van~den Broek, E., M.~J.~J. Dijkstra, O.~Krijgsman, D.~Sie, J.~C. Haan,
  J.~J.~H. Traets, M.~A. van~de Wiel, I.~D. Nagtegaal, C.~J.~A. Punt,
  B.~Carvalho, B.~Ylstra, S.~Abeln, G.~A. Meijer, and R.~J.~A. Fijneman (2015,
  09).
\newblock High prevalence and clinical relevance of genes affected by
  chromosomal breaks in colorectal cancer.
\newblock {\em PLOS ONE\/}~{\em 10\/}(9), 1--14.

\bibitem[\protect\citeauthoryear{Westfall and Wolfinger}{Westfall and
  Wolfinger}{1997}]{WestWolf1997}
Westfall, P. and R.~Wolfinger (1997).
\newblock Multiple tests with discrete distributions.
\newblock {\em The American Statistician\/}~{\em 51\/}(1), 3--8.

\end{thebibliography}


\appendix
\section{Additional procedures}

\subsection{A rescaled BH procedure}\label{sec:RBH}

The procedure  [RBH] (rescaled-BH) is defined as the step-up procedure using the critical values $\tau_k=\lambda_\alpha k/m$, $1\leq k \leq m$, where $\lambda_\alpha = \max\{ \lambda\in [0,1]\::\: \Psi(\lambda_\alpha)\leq \alpha\}$ for
$$
\Psi(\lambda)= \min\left(\lambda, \max_{1\leq k \leq m}\left( \frac{1}{k}  \sum_{i=1}^m \frac{F_i\left(\lambda k/m\right)}{1-F_i\left(\lambda\right)} \right)\right).
$$
The following result is straightforward from Theorem~\ref{th:step-updown} (SU part).

\begin{corollary} In the model of section~\ref{sec:model} with the additional assumption \eqref{equ-sousunif}, we have $\forall P \in\mathcal{P}$, $\FDR(\RBH,P)\leq \alpha$.
\end{corollary}

Moreover, if $\alpha$ is such that the equality $\Psi(\lambda_\alpha) = \alpha$ holds true, then $\lambda_\alpha\geq \Psi(\lambda_\alpha) = \alpha$ and [RBH] always dominates [BH] in terms of critical values and therefore rejects at least as many hypotheses. 

\subsection{ A discrete BR procedure}\label{sec:DBR}

The procedure [DBR-$\lambda$] (discrete BR) is defined as the step-up procedure  ${\SU}(\tau)$ using the critical values defined  in the following way $\tau_m =\max\left\{ t\in \mathcal{A}\::\: t \leq  \left((1-\lambda) m\alpha\right) \wedge \lambda \right\}$ (that is, the discretised last critical value of [BR-$\lambda$]) and 
\begin{align}
\tau_k &= \max\left\{ t\in \mathcal{A}\::\: t\leq \lambda,\:  \left(  F\left(t\right)\right)_{(1)}+\dots + \left(  F\left(t\right)\right)_{(m-k+1)}\leq \alpha k(1-\lambda)\right\},\:
1\leq k\leq m-1,\nonumber
\end{align}
where each $\left( F\left(t\right)\right)_{(j)}$ denotes the $j$-th largest elements of the set $\left\{ F_i\left(t\right) , 1\leq i \leq m \right\}$.
The following result is straightforward from Theorem~\ref{th:step-updown} (SU part).

\begin{corollary} In the model of section~\ref{sec:model} with the additional assumption \eqref{equ-sousunif}, we have $\forall P \in\mathcal{P}$, $\FDR(\DBR,P)\leq \alpha$. Moreover, the set of nulls rejected by [DBR-$\lambda$] is larger than the one of [BR-$\lambda$] (almost surely), with equality (almost surely) under \eqref{equ-unifgeneral} and $F_i=F_j$ for all $i\neq j$.
\end{corollary}

%

\section{Proof of Theorem~\ref{th:step-updown}}\label{sec:proofs}

\subsection{Lemmas for step-down and step-up procedures}

Let us introduce the following modifications of ${\SU}(\tau)$ :
\begin{itemize}
\item[\textbullet] ${\SU}'(\tau) = {\SU}(\tau') $ the step-up with $m$ critical values defined by $(\tau'_1,\dots,\tau'_m)=(\tau_2, \dots, \tau_m, \tau_m)$;
\item[\textbullet] for some given index $i\in\{1,\dots,m\}$, ${\SU}'^{-i}(\tau) = {\SU}(\tau'^{-i}) $ the step-up with $m-1$ critical values defined by $(\tau'^{-i}_1,\dots,\tau'^{-i}_{m-1})=(\tau_2, \dots, \tau_m)$ and restricted to the $p$-values of the set $\{p_{j},j\neq i\}$.
\end{itemize}

The following lemma holds (classical from \cite{FZ2006} and proved in Appendix~\ref{sec:supp}):

\begin{lemma}\label{lem:SU:proof}
For all  $i\in\{1,\dots,m\}$, the following assertions are equivalent: (i) $p_i\leq {\tau}_{\hat{k}}$; (ii) $p_i\leq {\tau}_{\hat{k}'^{-i}+1}$; (iii) $\wh{k}'^{-i}+1=\wh{k}$,
where $\wh{k}'^{-i}$ denotes the number of rejected hypotheses of the procedure ${\SU}'^{-i}(\tau)$. 
Moreover, we have $\{p_i > {\tau}_{m}\}\subset \{\wh{k}'=\wh{k}'^{-i}\}$, where $\wh{k}'$ denotes the number of rejected hypotheses of the procedure ${\SU}'(\tau)$. 
\end{lemma}

Let us introduce the following modifications of ${\SD}(\tau)$ :
\begin{itemize}
\item[\textbullet] for some given index $i\in\{1,\dots,m\}$, ${\SD}^{-i}(\tau) = {\SD}(\tau^{-i}) $ the step-down procedure with $m-1$ critical values defined by $(\tau^{-i}_1,\dots,\tau^{-i}_{m-1})=(\tau_1, \dots, \tau_{m-1})$ and restricted to the $p$-values of the set $\{p_{j},j\neq i\}$.
\item[\textbullet] for some given index $i\in\{1,\dots,m\}$, ${\SD}'^{-i}(\tau) = {\SD}(\tau'^{-i}) $ the step-down procedure with $m-1$ critical values defined by $(\tau'^{-i}_1,\dots,\tau'^{-i}_{m-1})=(\tau_2, \dots, \tau_m)$ and restricted to the $p$-values of the set $\{p_{j},j\neq i\}$.
\end{itemize}

The following lemma holds (classical from \cite{FZ2006} and proved in Appendix~\ref{sec:supp}):

\begin{lemma}\label{lemstepdown}
For all  $i\in\{1,\dots,m\}$, the following assertions are equivalent: (i) $p_i\leq {\tau}_{\tilde{k}}$; (ii) $p_i\leq {\tau}_{\tilde{k}+1}$; (iii) $p_i\leq {\tau}_{\tilde{k}^{-i}+1}$; (iv) $\tilde{k}'^{-i}+1=\tilde{k}$,
where $\wt{k}^{-i}$ is the number of rejections of ${\SD}^{-i}(\tau)$ and $\wt{k}'^{-i}$ is the number of rejections of ${\SD}'^{-i}(\tau)$.
Moreover, we have $\{p_i > {\tau}_{\tilde{k}^{-i}+1}\}\subset \{\wt{k}=\wt{k}^{-i}\}$.
\end{lemma}

\subsection{Proof of Theorem~\ref{th:step-updown}, step-up part}

By using Lemma~\ref{lem:SU:proof} $(ii)$ and $(iii)$ and independence, we easily obtain

\begin{align}
\FDR(\SU(\tau)) &= \sum_{i\in\cH_0}  \E\left( \frac{\ind{p_i\leq {\tau}_{\hat{k}}}}{\wh{k}}\right) 
=\sum_{i\in\cH_0}  \E\left( \frac{\ind{p_i\leq {\tau}_{\hat{k}'^{-i}+1}}}{\wh{k}'^{-i}+1}\right) 
\leq \sum_{i\in\cH_0}  \E\left(  \frac{F_i\left({\tau}_{\hat{k}'^{-i}+1}\right)}{\wh{k}'^{-i}+1}\right), \nonumber
\end{align}
where the last expectation is taken only with respect to $(p_j,j\neq i)$.
Now, on the one hand,
\begin{align*}
  \sum_{i\in\cH_0}  \E\left(  \frac{F_i\left({\tau}_{\hat{k}'^{-i}+1}\right)}{\wh{k}'^{-i}+1} \right)\leq 
 \sum_{i=1}^m  \E\left(  \frac{F_i\left({\tau}_{\hat{k}'^{-i}+1}\right)}{\wh{k}'^{-i}+1} \right)\leq \sum_{i=1}^m  \max_{1\leq k \leq m} \frac{F_i(\tau_k)}{k}.
\end{align*}
Next, on the other hand, by using again the independence,
\begin{align*}
 \sum_{i\in\cH_0}  \E\left( \frac{F_i\left({\tau}_{\hat{k}'^{-i}+1}\right)}{\wh{k}'^{-i}+1} \right) 
&\leq \sum_{i\in\cH_0}  \E\left( \frac{F_i\left({\tau}_{\hat{k}'^{-i}+1}\right)}{1-F_i\left({\tau}_{m}\right)} \frac{\ind{p_i> {\tau}_{m}}}{\wh{k}'^{-i}+1}\right) \\
&= \sum_{i\in\cH_0}  \E\left( \frac{F_i\left({\tau}_{\hat{k}'+1}\right)}{1-F_i\left({\tau}_{m}\right)} \frac{\ind{p_i> {\tau}_{m}}}{\wh{k}'+1} \ind{\wh{k}'+1\leq m}\right), 
\end{align*}
where the latter equality comes from the last assertion of Lemma~\ref{lem:SU:proof}.
Now, since ${\tau}_{\hat{k}'+1}\leq \tau_m$, we have that the last display is smaller or equal to 
\begin{align*}
&   \E\left( \sum_{i=1}^m  \frac{F_i\left({\tau}_{\hat{k}'+1}\right)}{1-F_i\left({\tau}_{m}\right)} \frac{\ind{p_i> {\tau}_{\hat{k}'+1}}}{\wh{k}'+1} \ind{\wh{k}'+1\leq m}\right) 
\leq  \max_{0\leq k \leq m-1}\max_{ \substack{ A \subset \{1,\dots,m\}\\ |A|=m-k }}  \sum_{i\in A}\frac{F_i\left({\tau}_{k+1}\right)}{1-F_i\left({\tau}_{m}\right)} \frac{1}{k+1} , 
\end{align*}
by taking the maximum over all the possible realizations of the set $A=\{1\leq i \leq m \::\: p_i> {\tau}_{\hat{k}'+1} \}= \{1\leq i \leq m \::\: p_i> {\tau'}_{\hat{k}'} \}$ which is the index set corresponding to the non-rejected null hypotheses of $\SU(\tau')$ (the latter being by definition of cardinality $m-\hat{k}'$).
This concludes the proof.

\subsection{Proof of Theorem~\ref{th:step-updown}, step-down part}\label{sec:proofsd}

It is similar to the step-up case, but relies now on Lemma~\ref{lemstepdown} $(iii)$ and $(iv)$ (and still independence): we have
\begin{align*}
\FDR(\SD(\tau)) = \sum_{i\in\cH_0}  \E\left( \frac{\ind{p_i\leq {\tau}_{\tilde{k}}}}{\wt{k}}\right) =\sum_{i\in\cH_0}  \E\left( \frac{\ind{p_i\leq {\tau}_{\tilde{k}^{-i}+1}}}{\tilde{k}'^{-i}+1}\right)
\leq \sum_{i\in\cH_0}  \E\left(  \frac{F_i\left({\tau}_{\tilde{k}^{-i}+1}\right)}{\tilde{k}^{-i}+1}\right),
\nonumber
\end{align*}
which gives the first part of the bound.
Next, by using independence and the last assertion of Lemma~\ref{lemstepdown}, we obtain
\begin{align*}
\sum_{i\in\cH_0}  \E\left(  \frac{F_i\left({\tau}_{\tilde{k}^{-i}+1}\right)}{\tilde{k}^{-i}+1}\right)&\leq \sum_{i\in\cH_0}  \E\left( \frac{F_i\left({\tau}_{\tilde{k}^{-i}+1}\right)}{1-F_i\left({\tau}_{\tilde{k}^{-i}+1}\right)} \frac{\ind{p_i>{\tau}_{\tilde{k}^{-i}+1}}}{\tilde{k}^{-i}+1}\right) \\
&\leq   \E\left( \sum_{i=1}^m\frac{F_i\left({\tau}_{\tilde{k}+1}\right)}{1-F_i\left({\tau}_{\tilde{k}+1}\right)} \frac{\ind{p_i>{\tau}_{\tilde{k}+1}}}{\tilde{k}+1} \ind{\tilde{k}+1\leq m}\right) \\
&\leq   \E\left(\max_{0\leq k \leq m-1}\max_{ \substack{A \subset \{1,\dots,m\}\\ |A|=m-k }}  \sum_{i\in A}\frac{F_i\left({\tau}_{k+1}\right)}{1-F_i\left({\tau}_{k+1}\right)} \frac{1}{k+1} \right), 
\end{align*}
 because $\{1\leq i \leq m \::\: p_i>{\tau}_{\tilde{k}+1}\}$ is equal to $\{1\leq i \leq m \::\: p_i>{\tau}_{\tilde{k}}\}$, which is the set of non-rejected hypotheses of $\SD(\tau)$. Since $\SD(\tau)$ rejects exactly $\tilde{k}$ hypotheses, the proof is completed. 

\section{Supplement}\label{sec:supp}

\subsection{Proofs for lemmas comparing procedures}\label{sec:proof:lemmaconserv}

The lemmas presented here rely on the fact that, there is almost surely no $p$-value in $[0,1]\backslash \mathcal{A}$ (both in the continuous and discrete cases). All  symbols ``$=$" or ``$\subset$" are intended to be valid almost surely in this section.

A result which will be extensively used in the proofs of this section is the following one : for $p$-values valued in the set $ \mathcal{A}$, then the step-up procedure with critical values $\tau_k$, $1\leq k \leq m$, has the same rejection set as the step-up procedure with critical values $\xi_k = \max \left\{t\in \mathcal{A}\::\: t \leq \tau_k\right\}$, $1\leq k \leq m$.
This fact comes from the simple following observation : for all $k$,
\begin{align*}
\{1\leq i \leq m\::\:  p_i\leq \tau_k\} &= \{1\leq i \leq m\::\:  p_i\in \mathcal{A},\:p_i\leq \tau_k\}\\
& = \{1\leq i \leq m\::\:  p_i\in \mathcal{A},\:p_i\leq \xi_k\} =  \{1\leq i \leq m\::\:  p_i\leq \xi_k\}.
\end{align*}
The $\xi_k$'s are called the ``effective" critical values of $\SD(\tau)$ or $\SU(\tau)$ in the sequel.

\subsubsection{Proof of Lemma~\ref{lem:HeyseBH}}

The  effective critical values of the BH procedure are given by $\xi_k = \max \left\{t\in \mathcal{A}\::\: t \leq \alpha k/m\right\}$, $1\leq k \leq m$. If \eqref{equ-sousunif} holds, then $\overline{F}(t)\leq t$ and each $\xi_k$ is clearly smaller than 
 the $k$-th critical values of [Heyse]. This implies that the rejection set of [Heyse] is larger than the one of [BH]. Conversely, under \eqref{equ-unifgeneral} and if $F_i=F_j=\overline{F}$ for all $i\neq j$, we always have $\overline{F}(t)=F_i(t)=t$ for $t\in\mathcal{A}$. This implies that the $\xi_k$'s are the critical values of [Heyse] and shows the reversed inclusion. 
 

\subsubsection{Proof of Lemmas~\ref{lem:DBH-SU-BH} and~\ref{lem:DBH-SD-BH}}\label{sec:prooflem:DBH-SU-BH}

Let $\tau_k$, $1\leq k \leq m$, be the critical values of [DBH-SU]. Let $\xi_k = \max \left\{t\in \mathcal{A}\::\: t \leq \frac{\alpha}{1+\alpha}\frac{k}{m} \right\}$ be the effective  critical values of the [BH] procedure at level $\alpha/(1+\alpha)$. 
Now, for all $t \in [0,1] $, we have by \eqref{equ-sousunif},
\begin{align}
\FbarSU(t) &= \frac{1}{m} \sum_{i=1}^m \frac{F_i(t)}{1-F_i(\tau_m)} \le 
\frac{t}{m} \sum_{i=1}^m \frac{1}{1-F_i(\tau_m)} =t \cdot (1+ \FbarSU(\tau_m)) \le t \cdot (1+\alpha), \label{eq:bound:fbarsu}
\end{align}
where the last inequality follows from the definition of $\tau_m$. Thus we have $\FbarSU(\xi_m)\le \alpha$, which in turn 
 implies $\xi_m \le \tau_m$. Additionally, the bound \eqref{eq:bound:fbarsu} yields for $1 \le k < m$
\begin{align*}
\tau_k &=  \max\left\{ t\in \mathcal{A}\::\: t\leq \tau_m,\: \FbarSU(t)  \leq \alpha k/m \right\}\\
&\geq  \max\left\{ t\in \mathcal{A}\::\: t\leq \tau_m,\: t (1+\alpha) \leq \alpha k/m \right\}\\
&=\max\left\{ t\in \mathcal{A}\::\:  t (1+\alpha) \leq \alpha k/m \right\}\\
&=\xi_k,
\end{align*}
where we used that $\xi_m \le \tau_m$. This proves Lemma~2.
The proof of Lemma~3 is analogue and is left to the reader.

\subsubsection{Proof of Lemma~\ref{lem:adapt}}

Let us first  focus on the case $(i)$ and denote by $\tau_k$, $1\leq k \leq m$, the critical values of [A-DBH-SU]. 
From \eqref{equ-sousunif}, we have for $1\leq k\leq m-1$,
\begin{align*}
\tau_k &\geq  \max\left\{ t\in \mathcal{A}\::\: t\leq \tau_m,\:   t\leq \alpha k(1-\tau_m) /(m-k+1) \right\}\\
&=  \max\left\{ t\in \mathcal{A}\::\: t\leq\left(  (1-\tau_m)\frac{ \alpha k}{m-k +1} \right)\wedge \tau_m\right\},
\end{align*}
which correspond to the effective critical values of [BR-$\lambda$] with $\lambda=\tau_m$.
Now consider the case $(ii)$ and denote again by $\tau_k$, $1\leq k \leq m$, the critical values of [A-DBH-SD]. From \eqref{equ-sousunif}, we have for $1\leq k\leq m$,
$$
\tau_k \geq \max\left\{ t\in \mathcal{A}\::\:   (m-k+1) t/(1-t) \leq \alpha k\right\}=  \max\left\{ t\in \mathcal{A}\::\:  t \leq \alpha k/(m-k(1-\alpha)+1)\right\}
$$
which correspond to the effective critical values of [GBS]. This implies the result.

\subsection{Proofs of technical lemmas for step-down and step-up procedures}


\subsubsection{Proof of Lemma~\ref{lem:SU:proof}}



First note 
that for any step-up procedure
 $$
 \wh{k}= \max\left\{k\in\{0,1,...,m\}\::\: \sum_{i=1}^m \ind{p_i\leq {\tau}_k} \geq k\right\},
 $$
 which is sometimes more handy, because this definition avoids to rely explicitly on the order statistics of the $p$-values.  

Now, it is not difficult to check that $\wh{k}'^{-i}\geq \wh{k}-1$ always holds:  this comes from the inequality
$$
\wh{k}-1=\sum_{j=1}^m \ind{p_j\leq \tau_{\hat{k}}} -1 \leq \sum_{j\neq i} \ind{p_j\leq \tau_{\hat{k}}} =  \sum_{j\neq i} \ind{p_j\leq \tau'^{-i}_{\hat{k}-1}},
$$
because $ \tau'^{-i}_{\ell-1}=\tau_{\ell}$ for $\ell \in\{2,\dots,m\}$ (note that we can assume without loss of generality $\wh{k}\geq 1$ here).
This means that $(i)$ implies $(ii)$. Now, when $p_i\leq {\tau}_{\hat{k}'^{-i}+1}$, we have
$$
\hat{k}'^{-i} =  \sum_{j\neq i} \ind{p_j\leq \tau'^{-i}_{\hat{k}'^{-i}}}   =\sum_{j\neq i} \ind{p_j\leq \tau_{\hat{k}'^{-i}+1}}  = \sum_{j=1}^m \ind{p_j\leq \tau_{\hat{k}'^{-i}+1}} -1
$$
which implies $\hat{k}'^{-i} +1 \leq \sum_{j=1}^m \ind{p_j\leq \tau_{\hat{k}'^{-i}+1}}$ and thus $\wh{k}'^{-i}+1\leq \wh{k}$. Since, again, $\wh{k}'^{-i}\geq \wh{k}-1$ always holds, we have
$\wh{k}'^{-i}+1 = \wh{k}$. Hence, $(ii)$ implies $(iii)$. Now, if  $\wh{k}'^{-i}+1 = \wh{k}$, we have
\begin{align*}
\ind{p_i\leq \tau_{\hat{k}}} &= \sum_{j=1}^m \ind{p_j\leq \tau_{\hat{k}}} - \sum_{j\neq i} \ind{p_j\leq \tau_{\hat{k}}}  = \wh{k} - 
 \sum_{j\neq i} \ind{p_j\leq \tau_{\hat{k}'^{-i}+1}} \\
 &= \wh{k} - 
 \sum_{j\neq i} \ind{p_j\leq \tau'^{-i}_{\hat{k}'^{-i}}} = \wh{k} -  \wh{k}'^{-i} =1,
 \end{align*}
by definition of $\tau'^{-i}$, which gives that $(iii)$ implies $(i)$. 
Now, to prove the last statement, we first note that $\wh{k}'\geq\wh{k}'^{-i}$ always holds. Furthermore, if $p_i > {\tau}_{m}$ let us prove $\wh{k}'\leq\wh{k}'^{-i}$. First, $\wh{k}'=m$ is impossible because $p_i$ is above ${\tau}_{m}$ and thus $p_i$ cannot be rejected by $\SU'(\tau)$. Hence, $\wh{k}'\leq m-1$ and thus $\tau'^{-i}_{\hat{k}'}$ is well defined. 
Now, since $p_i > {\tau}_{m}$, we obtain
$$
\sum_{j\neq i} \ind{p_j\leq \tau'^{-i}_{\hat{k}'}} = \sum_{j\neq i} \ind{p_j\leq \tau'_{\hat{k}'}}= \sum_{j=1}^m \ind{p_j\leq \tau'_{\hat{k}'}} = \wh{k}',
$$
which implies $\wh{k}'\leq\wh{k}'^{-i}$ by definition of ${\SU}'^{-i}(\tau)$.

\subsubsection{Proof of Lemma~\ref{lemstepdown}}

First note that 
for any step-down procedure
 $$
 \wt{k}= \max\left\{k\in\{0,1,...,m\}\::\: \forall k'\leq k, \: \sum_{i=1}^m \ind{p_i\leq {\tau}_{k'}} \geq k'\right\}.
 $$

Now, we check that $\wt{k}'^{-i}+1\geq \wt{k}$ always holds. Since $\sum_{j\neq i} \ind{p_j\leq \tau'^{-i}_{\tilde{k}'^{-i}+1}} < \tilde{k}'^{-i}+1$, we have 
$$
\sum_{j=1}^m \ind{p_j\leq \tau_{\tilde{k}'^{-i}+2}}  \leq 1+ \sum_{j\neq i} \ind{p_j\leq \tau'^{-i}_{\tilde{k}'^{-i}+1}} < \wt{k}'^{-i}+2,
$$
which gives $\wt{k}<\wt{k}'^{-i}+2$ by definition of $\wt{k}$ and thus $\wt{k}\leq \wt{k}'^{-i}+1$.
Next, if $p_i\leq {\tau}_{\tilde{k}}$, we have 
$$
\sum_{j\neq i} \ind{p_j\leq \tau'^{-i}_{\tilde{k}}} =\sum_{j\neq i} \ind{p_j\leq \tau_{\tilde{k}+1}} = \sum_{j=1}^m \ind{p_j\leq \tau_{\tilde{k}+1}} -1 < \wt{k}+1 -1,
$$ 
so that $\wt{k}>\wt{k}'^{-i}$ and thus $\wt{k}\geq \wt{k}'^{-i}+1$. This proves that $(i)$ implies $(iv)$. 
Next, if $p_i > {\tau}_{\tilde{k}^{-i}+1}$, then 
$$
\sum_{j=1}^m \ind{p_j\leq {\tau}_{\tilde{k}^{-i}+1}} = \sum_{j\neq i} \ind{p_j\leq {\tau}_{\tilde{k}^{-i}+1}} = \sum_{j\neq i} \ind{p_j\leq {\tau}^{-i}_{\tilde{k}^{-i}+1}} < \wt{k}^{-i}+1,
$$
which entails $\wt{k}<\wt{k}^{-i}+1$ and thus $\wt{k}\leq \wt{k}^{-i}$.
This proves $\wt{k}\neq \wt{k}'^{-i}+1$. Hence, $(iv)$ implies $(iii)$. 
The fact that  $(iii)$ implies $(ii)$ is obvious because $\wt{k}\geq \wt{k}^{-i}$ always holds.
Finally, we merely check that $\wt{k}$ is such that
$$
\wt{k} = \sum_{j=1}^m \ind{p_j\leq \tau_{\tilde{k}}} =  \sum_{j=1}^m \ind{p_j\leq \tau_{\tilde{k}+1}} ,
$$
which means that the set of $p$-values rejected at threshold $\tau_{\tilde{k}}$ is the same as the set of $p$-values rejected at threshold $\tau_{\tilde{k}+1}$. This gives that $(ii)$ implies $(i)$.
For the last assertion, it has been proved in the above reasoning while showing that $(iv)$ implies $(iii)$.

\subsection{Table for the simulations}

\begin{table}
\caption{Average power of FDR procedures for $N=25$ (see Section~\ref{sec:simu}). 	\label{tab:Appendix1}}
	 \centering
	\small
		\begin{tabular}{cccc||ccccc}
$m $    & $m_3 $   & $m_1$    & $q$    & [BH] & [Heyse] & [DBH-SU] & [A-DBH-SU] & [A-DBH-SD] \\
\midrule
		800   & 80    & 144   & 0.15  & 0.0000 & 0.0004 & 0.0003 & 0.0003 & 0.0004 \\
	&       & 144   & 0.25  & 0.0004 & 0.0197 & 0.0177 & 0.0177 & 0.0135 \\
	&       & 144   & 0.4   & 0.0803 & 0.4425 & 0.4247 & 0.4247 & 0.4130 \\
	&       & 360   & 0.15  & 0.0000 & 0.0007 & 0.0006 & 0.0006 & 0.0007 \\
	&       & 360   & 0.25  & 0.0004 & 0.0244 & 0.0209 & 0.0209 & 0.0153 \\
	&       & 360   & 0.4   & 0.0803 & 0.4529 & 0.4509 & 0.4509 & 0.4487 \\
	&       & 576   & 0.15  & 0.0000 & 0.0009 & 0.0007 & 0.0007 & 0.0008 \\
	&       & 576   & 0.25  & 0.0004 & 0.0343 & 0.0259 & 0.0259 & 0.0231 \\
	&       & 576   & 0.4   & 0.0803 & 0.5367 & 0.4741 & 0.4741 & 0.4999 \\
	& 240   & 112   & 0.15  & 0.0000 & 0.0003 & 0.0003 & 0.0003 & 0.0002 \\
	&       & 112   & 0.25  & 0.0005 & 0.0276 & 0.0249 & 0.0249 & 0.0157 \\
	&       & 112   & 0.4   & 0.2148 & 0.5365 & 0.5012 & 0.5012 & 0.4951 \\
	&       & 280   & 0.15  & 0.0000 & 0.0003 & 0.0003 & 0.0003 & 0.0002 \\
	&       & 280   & 0.25  & 0.0005 & 0.0315 & 0.0272 & 0.0272 & 0.0175 \\
	&       & 280   & 0.4   & 0.2147 & 0.5758 & 0.5536 & 0.5536 & 0.5495 \\
	&       & 448   & 0.15  & 0.0000 & 0.0005 & 0.0003 & 0.0003 & 0.0004 \\
	&       & 448   & 0.25  & 0.0005 & 0.0372 & 0.0308 & 0.0308 & 0.0207 \\
	&       & 448   & 0.4   & 0.2145 & 0.5920 & 0.5741 & 0.5741 & 0.5775 \\
	& 640   & 32    & 0.15  & 0.0000 & 0.0002 & 0.0002 & 0.0002 & 0.0002 \\
	&       & 32    & 0.25  & 0.0010 & 0.0378 & 0.0341 & 0.0341 & 0.0174 \\
	&       & 32    & 0.4   & 0.4243 & 0.6174 & 0.5955 & 0.6828 & 0.6621 \\
	&       & 80    & 0.15  & 0.0000 & 0.0002 & 0.0002 & 0.0002 & 0.0002 \\
	&       & 80    & 0.25  & 0.0010 & 0.0388 & 0.0347 & 0.0347 & 0.0179 \\
	&       & 80    & 0.4   & 0.4242 & 0.6282 & 0.6128 & 0.6841 & 0.6638 \\
	&       & 128   & 0.15  & 0.0000 & 0.0002 & 0.0002 & 0.0002 & 0.0002 \\
	&       & 128   & 0.25  & 0.0010 & 0.0400 & 0.0354 & 0.0354 & 0.0183 \\
	&       & 128   & 0.4   & 0.4240 & 0.6353 & 0.6265 & 0.6854 & 0.6656 \\
	2000  & 200   & 360   & 0.15  & 0.0000 & 0.0002 & 0.0002 & 0.0002 & 0.0002 \\
	&       & 360   & 0.25  & 0.0001 & 0.0156 & 0.0142 & 0.0142 & 0.0100 \\
	&       & 360   & 0.4   & 0.0730 & 0.4486 & 0.4317 & 0.4317 & 0.4197 \\
	&       & 900   & 0.15  & 0.0000 & 0.0002 & 0.0002 & 0.0002 & 0.0002 \\
	&       & 900   & 0.25  & 0.0001 & 0.0192 & 0.0166 & 0.0166 & 0.0125 \\
	&       & 900   & 0.4   & 0.0730 & 0.4517 & 0.4511 & 0.4511 & 0.4509 \\
	&       & 1440  & 0.15  & 0.0000 & 0.0003 & 0.0002 & 0.0002 & 0.0002 \\
	&       & 1440  & 0.25  & 0.0001 & 0.0286 & 0.0211 & 0.0211 & 0.0165 \\
	&       & 1440  & 0.4   & 0.0730 & 0.5402 & 0.4684 & 0.4684 & 0.4984 \\
	& 600   & 280   & 0.15  & 0.0000 & 0.0002 & 0.0002 & 0.0002 & 0.0002 \\
	&       & 280   & 0.25  & 0.0001 & 0.0239 & 0.0213 & 0.0213 & 0.0115 \\
	&       & 280   & 0.4   & 0.2058 & 0.5350 & 0.4988 & 0.4988 & 0.4960 \\
	&       & 700   & 0.15  & 0.0000 & 0.0002 & 0.0002 & 0.0002 & 0.0002 \\
	&       & 700   & 0.25  & 0.0001 & 0.0290 & 0.0239 & 0.0239 & 0.0132 \\
	&       & 700   & 0.4   & 0.2058 & 0.5750 & 0.5590 & 0.5590 & 0.5516 \\
	&       & 1120  & 0.15  & 0.0000 & 0.0002 & 0.0002 & 0.0002 & 0.0002 \\
	&       & 1120  & 0.25  & 0.0001 & 0.0350 & 0.0283 & 0.0283 & 0.0157 \\
	&       & 1120  & 0.4   & 0.2057 & 0.5908 & 0.5739 & 0.5739 & 0.5761 \\
	& 1600  & 80    & 0.15  & 0.0000 & 0.0001 & 0.0001 & 0.0001 & 0.0000 \\
	&       & 80    & 0.25  & 0.0003 & 0.0379 & 0.0342 & 0.0342 & 0.0126 \\
	&       & 80    & 0.4   & 0.4223 & 0.6196 & 0.5928 & 0.6860 & 0.6591 \\
	&       & 200   & 0.15  & 0.0000 & 0.0001 & 0.0001 & 0.0001 & 0.0000 \\
	&       & 200   & 0.25  & 0.0003 & 0.0387 & 0.0351 & 0.0351 & 0.0131 \\
	&       & 200   & 0.4   & 0.4222 & 0.6281 & 0.6152 & 0.6869 & 0.6602 \\
	&       & 320   & 0.15  & 0.0000 & 0.0001 & 0.0001 & 0.0001 & 0.0000 \\
	&       & 320   & 0.25  & 0.0003 & 0.0396 & 0.0360 & 0.0360 & 0.0137 \\
	&       & 320   & 0.4   & 0.4220 & 0.6327 & 0.6278 & 0.6877 & 0.6617 \\
	\bottomrule
	\end{tabular}%
\end{table}%

\end{document}